\documentclass[leqno,12pt]{amsart}
\usepackage{amscd,amsfonts,amsmath,amssymb,amstext,amsthm}
\newtheorem{theorem}{Theorem}[section]
\newtheorem{lemma}[theorem]{Lemma}

\begin{document}
\title{On Common Zeros
of
Eigenfunctions \\ of the Laplace Operator}
\author{Dmitri Akhiezer and Boris Kazarnovskii}
\address {Institute for Information Transmission Problems \newline
19 B.Karetny per.,127994, Moscow, Russia,\newline
{\rm D.A.:} 
{\it akhiezer@iitp.ru},  
{\rm B.K.:} 
{\it kazbori@gmail.com}.}
\subjclass{53C30, 58J05}
\keywords{Homogeneous Riemannian manifold, Laplace operator,  
equivariant map, kinematic formula}
\thanks{Research supported by Russian Foundation of Sciences, 
project No. 14-50-00150}
\begin{abstract}
We consider the eigenfunctions of the Laplace operator $\Delta $ on a compact
Riemannian manifold $M$ of dimension $n$. For $M$ homogeneous with irreducible
isotropy representation and for a fixed eigenvalue $\lambda $ of $\Delta $
we find the average number of common zeros of $n$ eigenfunctions. 
It turns out that, up to a constant depending on $n$, this number equals
$\lambda^{n/2}{\rm vol}\,M$, the expression known from the celebrated Weyl's
law.
To prove this we
compute the volume of the image of $M$ under an
equivariant immersion into a sphere.
\end{abstract}
\renewcommand{\subjclassname}
{\textup{2010} Mathematics Subject Classification}
\maketitle
\renewcommand{\thefootnote}{}
\maketitle

\section{Introduction}\label{intro}
Let $M$ be a compact Riemannian manifold 
without boundary, $n = {\rm dim}\, M$. Fix an eigenvalue
$\lambda $ of the Laplace operator $\Delta $ on $M$, such that the
corresponding eigenspace
$$W_\lambda =\{ u\in C^\infty (M, \mathbb R)\ \vert \ \Delta u + \lambda u = 0\}$$
has sufficiently large dimension.
Moreover, assume that for some eigenfunctions $u_1,\ldots,u_n \in W_\lambda$ 
the set of common zeros 
$$Z(u_1,\ldots,u_n) =\{x\in M \ \vert \ u_1(x) = \ldots = u_n(x) = 0\}$$
is finite.
As we will see later, this happens quite often for homogeneous spaces
of compact Lie groups.
Our goal is to evaluate the number of points in $Z(u_1,\ldots,u_n)$.

The zero set of an eigenfunction is called a nodal set. The connected
components of its
complement are called nodal domains. 
According to the classical Courant's theorem \cite{CD},
the number of nodal domains determined 
by $k$-th eigenfunction is
at most $k$, see also \cite{Ch}, \cite{Pl}. 
For the standard sphere $S^2$ the eigenfunctions are spherical harmonics.
A spherical harmonic of degree $m$ has  
eigenvalue $\lambda _ m = m(m+1)$ of multiplicity $2m+1$, so the
number of nodal domains is smaller than or equal to $k=m^2 + 1$.
Therefore the nodal set has at most $m^2$ connected components.
On the other hand, it follows from B\'ezout's theorem that
the number of points in $Z(u_1,u_2)$ for two generic spherical harmonics
of degree $m$ does not exceed $2m^2$, see Sect.~\ref{bez}. 
Now, if $u_1$ is fixed, one can define $u_2$ as a small perturbation
of $u_1$
by a rotation. Let $p$ be the number of connected
components of the nodal set $u_1 = 0$.
Take one of them and observe that it 
intersects some connected component of the nodal set $u_2 = 0$.
Moreover, by
Jordan theorem applied to a loop of the former connected
component it follows that the number
of intersection points is at least 2. 
Altogether, we get $2p$
points in the zero set $Z(u_1,u_2)$, proving the Courant's estimate $p\le m^2$. 

Motivated by this example,
one is led to the problem of estimating $\# Z (u_1,\ldots,u_n)$
from above. Such an estimate is a simple result for the sphere,
see Theorem~\ref{est}. The proof is based on the relation
between spherical harmonics and homogeneous harmonic polynomials
on the ambient vector space, so there is no 
generalization to arbitrary $M$.
Alternatively, for $M$ homogeneous and isotropy irreducible we  
find the average number  
of common zeros of eigenfunctions.
In \cite{Ar}, V.Arnold suggested to study the topology of the set of
common zeros for $m\le n$ eigenfunctions, see Problem 2003-10, p.\,174.
Under our assumptions, we consider this problem for $m=n$.

We now state our main result.
Let $K$ be a connected compact Lie group,
$L\subset K$ a closed subgroup, and $M = K/L$. 
Fix a $K$-invariant Riemannian metric $g$ on $M$. 
In what follows, we consider real valued functions on $M$,
for which the scalar product is given by 
$$(f_1, f_2) = \int _Mf_1(x)f_2(x)dx,$$
where $dx$ is the Riemannian measure.
Since the Laplace operator
$\Delta $ on $M$ is $K$-invariant, 
we have an orthogonal representation of $K$ in 
the eigenspace $W_\lambda $. From now on we tacitly assume that $\lambda >0$.
Let $\{0\}\ne H\subset W_{\lambda } $ be a
$K$-invariant subspace of dimension
$N $.

For $x\in M$ we have the linear functional $\alpha _x \in H^*$, defined
by $\alpha _x(u) = u(x)$, where $u\in H$.  
Take
an orthonormal basis $\{f_i\}_{i=1}^N$ of $H$
and identify $H^*$ with ${\mathbb R}^N$ via
$$(\mu _1, \ldots, \mu_N) \mapsto \mu_1f_1^* + \ldots + \mu_Nf_N^*,$$
where $f_i^*(f_j) = \delta _{ij}$. Then
$$\alpha _x = \sum f_i(x)f^*_i$$ and so the map $M \to H^*,\ x \mapsto 
\alpha _x,$ is written as
$$f=(f_1,\ldots,f_N): M \to {\mathbb R}^N.$$
Clearly, $f$ is equivariant with respect to the given action of $K$ on $M$
and the linear representation of $K$ in $H^* $ identified with ${\mathbb R}^N$.

We now restrict our attention to the class
of isotropy irreducible homogeneous spaces. 
This means that
the the subgroup $L\subset K$ acts
in the tangent space to $M$
via an irreducible representation. Isotropy irreducible
homogeneous spaces are listed in \cite{Ma} and \cite{Wo}. 
We note
that symmetric spaces of simple groups belong to this class.

In what follows, we denote by $\sigma _n$ the $n$-dimensional volume
of the unit sphere $S^n \subset \mathbb R^{n+1}$.              
We write $Z(U)$ for the set of common zeros of all functions from an 
$n$-dimensional subspace $U \subset H$. If $u_1, \ldots, u_n$
is a basis of $U$ then $Z(U) = Z(u_1, \ldots, u_n)$.
The average 
${\mathfrak M}\, \{\#Z(U)\}$
is taken over all $n$-dimensional subspaces $U$ of $H$, i.e., over
the Grassmanian ${\rm Gr}_n(H)$, with respect to the action of ${\rm SO}(N,
\mathbb R)$, endowed with the Haar measure of volume 1, on the set of these 
subspaces. 
\begin{theorem}\label{main}
Let $M = K/L$ be isotropy irreducible and 
let $H \subset W_\lambda$ be a $K$-invariant subspace. 
Then
$${\mathfrak M}\, \{\#Z(U)\} = \frac{2}{\sigma _n}
\Bigl(\frac {\lambda}{n}\Bigr)^{n/2}{\rm vol}\, M = 
c_n\lambda^{n/2}{\rm vol}\, M,$$
where $c_n$ depends only on the dimension of $M$.
In particular, for $M = S^n$ and for the 
eigenvalue $\lambda = m(m+n-1)$
one has
$${\mathfrak M}\, \{\#Z(U)\}
= 2 
\Bigl(\frac {m(m+n-1)}{n}\Bigr)^{n/2}.$$
\end{theorem}
\medskip
\noindent
We remark that in the latter case
the representation in $W_\lambda $ is irreducible, i.e., $H = W_\lambda$.

\medskip \noindent {\sl
The authors would like to thank S.Alesker for pointing out a 
reference for a relevant Crofton's formula. }

\section{Equivariant mappings and induced metrics }\label{}

We use the notations introduced in Sect.~\ref{intro}.
The scalar product $g(\xi,\eta)$ on a tangent space 
$T_x(M)$ defines a scalar product on the dual space $T_x^*(M)$,
denoted again by $g$. Furthermore,
for tangent vectors and covectors 
we write $(\xi,\eta) $ instead of 
$g(\xi,\eta)$ and $\vert \vert \xi \vert \vert ^2$ instead of $g(\xi,\xi)$.

\begin{theorem} \label{diff} One has
$$\sum _{i=1}^N \, \vert \vert df_i \vert \vert^2 
= \frac{\lambda N}{{\rm vol}M}$$
everywhere on $M$.
\end{theorem} 
\begin{proof} Fix $k \in K$ and put $f_i^*(x) = f_i(k^{-1}x)$.
Then 
$$f_i^*(x) = \sum _{j}\, a_{ij}f_j(x),\ \ df_i^*(x) = \sum _{j}\, 
a_{ij}df_j(x),$$
where $(a_{ij})$ is an orthogonal matrix. Therefore
$$
\vert \vert df_i^*(x)\vert \vert ^2 
= \sum_{p,q}\,a_{ip}a_{iq}\cdot(df_p(x),df_q(x)),
$$
hence
$$
\sum _i\vert\vert df_i(x)\vert \vert ^2 
= \sum _i\vert \vert df_i^*(x)
\vert\vert^2 .
$$
Now,
for $\xi \in T_xM$ and $\eta = (dk)\xi \in T_{kx}M$
one has
$$df_i^*(kx)(\eta) = df_i(x)(\xi),$$
where $\xi \in T_xM$ and $\eta = (dk)\xi \in T_{kx}$.
Taking the maximum over the unit sphere in $T_xM$,
we get
$$\vert \vert df_i^*(kx)\vert \vert ^2  =
\vert \vert df_i(x)\vert \vert ^2.$$
Since $k$ is arbitrary, it follows that
$\sum _i \vert\vert df_i(x)\vert \vert ^2  
=D, $
where $D$ is a constant. 
To find $D$ write
$${\rm div}(f_i{\rm grad}f_i) = f_i\Delta f_i + \vert \vert {\rm grad}f_i
\vert \vert ^2 = - \lambda f_i^2 +\vert \vert df_i\vert \vert^2.$$
Summing up and integrating over $M$, we obtain
$$\lambda N =\lambda \sum_i\int_M\,f_i^2 dx = D\cdot{\rm vol}\,M ,$$
hence
$$D = \frac {\lambda N}{{\rm vol}\,M}.$$
\end{proof}

\begin{lemma}\label{forms} Let $l_k,\ k=1,\ldots,N$ be linear forms in $n$
variables $t_j$, such that $\sum l_k^2 = C\cdot \sum t_j^2$. Then
$\sum \vert \vert l_k\vert \vert ^2 = Cn$,
where $\vert \vert l \vert \vert := \sqrt{\sum b_j^2}$ for $l = \sum b_jt_j$.

\end{lemma}
\begin{proof} Write $l_k = \sum\, b_{kj}t_j$. Then 
$\sum_{k}\, b_{ki}b_{kj} = C\cdot \delta _{ij} $ by assumption. In particular,
$$\sum_{k=1}^N\, \vert \vert l_k \vert \vert ^2 = \sum_k\,\sum_j\, b_{kj}^2=
\sum_j\,\sum_k\, b_{kj}^2 = Cn.$$  
\end{proof}
\medskip 

For the convenience of the reader we give the proofs of some
known results.
In quantum mechanics, the next theorem for spherical harmonics
is called 
Uns\"old's theorem (1927).
The general case is also found in the literature,
see \cite{GHL}, Exercise 5.25 c),i), p. 261 and
p. 303.   
\begin{theorem}\label{rad}
One has
$$
\sum_{i=1}^N \,f_i^2 = R^2,
$$
where 
$R = \sqrt \frac{N}{{\rm vol}\, M},$ 
i.e., $f(M)$ is contained in the sphere $S^{N-1}(R)\subset {\mathbb R}^N$ 
of radius  
$R$. 
\end{theorem}
\begin{proof} 
Define $f^*_i(x)$ as in the proof of Theorem \ref{diff}.
The same argument using the orthogonality of $(a_{ij})$ shows
that
$$\sum _i\, f_i(x)^2 
= \sum _i \,f_i^*(x)^2.$$
By transitivity of $K$ on $M$ it follows that 
$$
\sum _{i=1}^N\, f_i^2 = R^2,
$$
where $R^2$ is a constant. Integrating over $M$ yields
$R^2\cdot {\rm vol}\, M = N$. 
\end{proof}
\bigskip
\noindent

Let $\mathfrak k$ and $\mathfrak l$ be the Lie algebras of $K$ and,
respectively,
$L$. If $M=K/L$ is isotropy irreducible, then the $L$-module
${\mathfrak k}/{\mathfrak l}$, defined by the adjoint action
of $L$, is irreducible. In particular, an $L$-invariant
intermediate Lie subalgebra between $\mathfrak l$ and $\mathfrak k$ coincides
with one of these two algebras. Therefore, for a subgroup $L_1\subset K$
containing $L$, we have one of the two
possibilities: either $L\subset L_1$ is a finite extension, or
$L_1 =K$.   

The first assertion in the following theorem is known,
see \cite{GHL}, Exercise 5.25 c), ii), p. 261 and p. 303.
\begin{theorem}\label{dil}
Assume that the isotropy representation of $L$
in the tangent space to $M$
is irreducible. Then:

\noindent 1) the mapping $f:M \to f(M)$ is a covering
of some degree $d$ and a local isometry up to a dilation; 

\noindent 2) the inverse image of Euclidean metric on $\mathbb R^N$
is given by
$$f^*(\sum \, dt_i^2) = \frac{\lambda N}{n{\rm vol}\,M} \cdot g$$ 
and the volume of $f(M) \subset S^{N-1}(R)$ is equal to
$${\rm vol}\,f(M) = \frac{1}{d}\cdot 
\Bigl(\frac{\lambda N}
{n {\rm vol}\,M} \Bigr)^{n/2}\cdot
{\rm vol}\, M.$$  
\end{theorem}
\begin{proof}  
1) Recall that $f:M \to {\mathbb R}^N$ is an equivariant map with respect
to the action of $K$ on the homogeneous space $M$
and a linear representation of $K$ in ${\mathbb R}^N$. Therefore $f(M)$
is one orbit of $K$ in ${\mathbb R}^N$. 
In particular, $f(M)$ is a manifold and a homogeneous space.

By our assumption,
an equivariant mapping from $M$ to any homogeneous space of $K$
is either a covering map or the map to one point.
But $\lambda > 0$, the functions $f_i$ are non-constant,
and so $f$ is a covering map.
By the irreducibility of the isotropy representation
it follows that there is only one up to a scalar factor
$K$-invariant metric on $M$. Hence
$$f^*(\sum\, dt_i^2) = C\cdot g.$$

\noindent 2) By Lemma \ref{forms} and Theorem \ref{diff} we have 
$$Cn = \sum \vert \vert df_i \vert \vert ^2 = \frac{\lambda N}{{\rm vol}\, M}.$$
The measure on $M$ induced by the metric $C\cdot g$
is $C^{n/2}\cdot dx$, and so we obtain the expression for ${\rm vol}\, f(M)$.

\end{proof}

\section{The kinematic formula}\label{cr}
The kinematic formula is proven by L.A.Santal\'o in his
book \cite{Sa}, see Sect.15.2 and 18.6 therein.
Another proof is given by R.Howard, see \cite{Ho}, Sect. 3.12.
The formula is valid for any space of constant curvature,
but we will use it only for the sphere. Also, we will need only the special 
case when two submanifolds of the sphere have complementary dimensions.
In geometric probability, the results of this type are called Crofton formulae.
Let $M$ and $L$ be two submanifolds of the sphere $S^{N-1} \subset 
{\mathbb R}^N$, $n= {\rm dim}\, M$, $l={\rm dim}\, L$, and $n+l = N-1$.
For the applications it is convenient to consider the sphere of radius $R$.
The kinematic density is given by the Haar measure $dg$ on 
${\rm SO}(N,\mathbb R)$.
The kinematic formula is the equality
$$\int _{M \cap g\cdot L \ne \emptyset}\#(M \cap g\cdot L)\cdot dg = 
C\, ({\rm vol}\, M)  \,
({\rm vol}\, L),$$
where $C$ is a constant independent of $M$ and $L$.
The constant is easy to find. Namely, draw two planes of dimensions 
$n+1$ and $l+1$ through the origin and take
for $M$ and $L$ the plane
sections. They are isometric to spheres of dimensions $n$ and $l$,
and the number of intersection points is 2 almost everywhere.
Therefore
$$C = \frac{2}{\sigma _n\sigma _lR^{N-1}}.$$
We will apply the kinematic formula for arbitrary $M$ and for a plane
section $L$ of complementary dimension through the origin. 
The kinematic formula means that the
average number of intersection points
of $M$ by such a
plane is equal to
$$\int _{M \cap g\cdot L \ne \emptyset}\#(M \cap g\cdot L)\cdot dg = 
\frac{2\,{\rm vol}\,M}{\sigma _nR^n}.$$

\medskip
\noindent
{\it Proof of Theorem ~\ref{main}}.
We have the covering $f:M \to f(M) \subset S^{N-1}$
of degree $d$. The set of common zeros for an $n$-dimensional
subspace $U \subset H \subset
W_\lambda $
can be written as
$$Z(U) = f^{-1}(f(M) \cap L)$$
for some $L$. If this set is finite then
$$\# Z(U) = d\cdot \#(f(M)\cap L).$$
For the average number of zeros the kinematic formula gives
$${\mathfrak M}\{\#Z(U)\} =
d\cdot \int_{f(M) \cap g\cdot L \ne 0} \#(f(M)\cap g\cdot L)dg =$$

$$= \frac{2d\cdot {\rm vol}\,f(M)}{\sigma _n R^n} = 
\frac{2}{\sigma _n}\Bigl(\frac{\lambda}{n}\Bigr)^{n/2}\cdot {\rm vol}\, M,$$
where the last equality follows from Theorems ~\ref{rad} and ~\ref{dil}, 2).

\hfill{$\square$}
 
\section{Upper bound in the case $M=S^n$}\label{bez}
Let $M = S^n$ be the unit sphere in ${\mathbb R}^{n+1}$.
An eigenfunction of
the Laplace operator on $M$ is the restriction   
of a homogeneous harmonic polynomial $p(t_1, \ldots, t_{n+1})$.
The polynomial is uniquely defined 
by the maximum principle for harmonic functions. If 
$m={\rm deg}\,p$
then $u = p\vert_{S^n}$ satisfies $\Delta u + \lambda u = 0$,
where $\lambda = m(m+n-1)$.
Let $u_i = p_i\vert_{S^n},\ {\rm deg}\,p_i = m_i,\ i=1,2,\dots,n$.
We call $u_1,\ldots,u_n$ generic if $p_1, 
\ldots, p_n$ have finitely many common zeros on
the complex affine hypersurface $\sum z_i^2 = 1$ in $\mathbb C^{n+1}$.
 
\begin{theorem}\label{est}
Let $u_1,\ldots,u_n$ be a generic system of spherical harmonics.
Then
$$\#Z(u_1,\ldots,u_n) \le 2\cdot m_1\cdot \ldots \cdot m_n.$$
\end{theorem}

\begin{proof} Consider 
the imbedding $\mathbb C^{n+1} \subset \mathbb P^{n+1}
(\mathbb C)$ and denote by
$X$ the intersection of hypersurfaces
$$p_1(t_1,\ldots ,t_{n+1}) = 0,\, \ldots ,\, p_n (t_1,\ldots, t_{n+1})= 0,\
t_1^2 +\ldots + t_{n+1}^2 = s^2$$
in $\mathbb P^{n+1}(\mathbb C)$ with homogeneous coordinates 
$(t_1:\ldots :t_{n+1}:s)$.
Then we have the inclusion of finite sets
$$X(\mathbb R) = Z(u_1,\dots,u_n) \subset X\cap {\mathbb C}^{n+1}.$$
It follows that the number of irreducible components 
of $X$ is greater than or equal to $\#Z(u_1,\ldots, u_n)$.
On the other hand, by B\'ezout's theorem this number does not
exceed the product of degrees $2\cdot \prod m_i$, 
see ~\cite{Fu}, Example 8.4.6.
\end{proof}

\section{Concluding remarks}
1.\ Let $\lambda _1, \lambda _2, \ldots, \lambda_n$ 
be $n$ not necessarily distinct
positive eigenvalues of $\Delta$. For each $\lambda _i$ choose a subspace $H_i
\subset W_{\lambda _i}$ as in Sect.~\ref{intro}, put $N_i = {\rm dim}\, 
H_i$, and consider
the equivariant map $f^{(i)} : M \to S^{N_i-1}$. Then the map
$$f: M \to S^{N_1-1} \times \ldots \times S^{N_n-1},\ 
\ \ x\mapsto(f^{(1)}(x), \ldots , f^{(n)}(x))$$
is equivariant with respect to the direct sum of representations
of $K$.
For $u_1 \in H_1, \ldots ,u_n \in H_n$ one can
mimick our definition in Sect.~\ref{intro} to get the average number
of zeros 
${\mathfrak M}\, \{\#Z(u_1, \ldots ,u_n)\}$.                
This amounts to averaging the number
of intersection points of $f(M)$ with the products $L_1 \times \ldots \times
L_n$, where $L_i$ is a hyperplane section of $S^{N_i-1}$ through the origin.
It is then plausible to conjecture that 
${\mathfrak M}\, \{\#Z(u_1,\ldots, u_n)\} $ equals
$$ \frac{2
\sqrt{\lambda_1\ldots\lambda_n}}{\sigma_n\,n^{n/2}}\cdot {\rm vol}\, M,$$
generalizing Theorem ~\ref{main}.
Unfortunately,
since there is no kinematic formula 
for the product of spheres, the proof in Sect. ~\ref{cr} does not work.

2.\ It would be interesting 
to have a lower bound for $\#Z(u_1, \ldots, u_n)$
at least for $M = S^n$. For $M=S^2$ one can easily construct two spherical
harmonics of degree $m$, such that $\#Z(u_1, u_2) = 2m$. Namely,
let $v_m$ be the so called zonal spherical harmonic, i.e.,
the spherical harmonic invariant under the group
of rotations around a given axis. The zeros of $v_m$
are located on $m$ circles orthogonal to the rotation axis.
Take another axis at angular distance $\alpha $ 
from the given one and denote by $v_m^\alpha $ the associate zonal spherical
harmonic of degree $m$.
The zeros of $v_m^{\alpha}$ are located on $m$ circles,
orthogonal to the new rotation axis.
If $\alpha $ is small enough then each circle in $v_m^{\alpha } = 0$
intersects exactly one circle in ${v_m = 0}$ in two points. Thus
$\#Z(v_m, v_m^{\alpha}) = 2m$.

\end {document}